\input amstex
\documentstyle{amsppt}
%
\catcode`@=11
\redefine\output@{%
  \def\break{\penalty-\@M}\let\par\endgraf
  \ifodd\pageno\global\hoffset=105pt\else\global\hoffset=8pt\fi  
  \shipout\vbox{%
    \ifplain@
      \let\makeheadline\relax \let\makefootline\relax
    \else
      \iffirstpage@ \global\firstpage@false
        \let\rightheadline\frheadline
        \let\leftheadline\flheadline
      \else
        \ifrunheads@ 
        \else \let\makeheadline\relax
        \fi
      \fi
    \fi
    \makeheadline \pagebody \makefootline}%
  \advancepageno \ifnum\outputpenalty>-\@MM\else\dosupereject\fi
}
\def\Beta{\mathchar"0\hexnumber@\rmfam 42}
\catcode`\@=\active
\nopagenumbers
\chardef\textvolna='176
\def\negskp{\hskip -2pt}

\chardef\degree="5E
\def\blue#1{#1}

\catcode`#=11\def\diez{#}\catcode`#=6
\catcode`&=11\catcode`&=4
\catcode`_=11\def\podcherkivanie{_}\catcode`_=8
\catcode`~=11\def\volna{~}\catcode`~=\active
\def\mycite#1{\cite{\blue{#1}}\immediate\special{ps:
     ShrHPSdict begin /ShrBORDERthickness 0 def}}
\def\myciterange#1#2#3#4{\cite{\blue{#2#3#4}}\immediate\special{ps:
     ShrHPSdict begin /ShrBORDERthickness 0 def}}
\def\mytag#1{%
    \tag#1}
\def\mythetag#1{\thetag{\blue{#1}}\immediate\special{ps:
     ShrHPSdict begin /ShrBORDERthickness 0 def}}
\def\myrefno#1{\no#1}
\def\myhref#1#2{\blue{#2}\immediate\special{ps:
     ShrHPSdict begin /ShrBORDERthickness 0 def}}
\def\myEarXivlink{\myhref{http://arXiv.org}{http:/\negskp/arXiv.org}}

\def\mytheorem#1{\csname proclaim\endcsname{Theorem #1}}
\def\mytheoremwithtitle#1#2{\csname proclaim\endcsname{Theorem #1#2}}

\def\mylemma#1{\csname proclaim\endcsname{Lemma #1}}
\def\mylemmawithtitle#1#2{\csname proclaim\endcsname{Lemma #1#2}}
\def\mythelemma#1{\blue{#1}\immediate\special{ps:
     ShrHPSdict begin /ShrBORDERthickness 0 def}}
\def\mycorollary#1{\csname proclaim\endcsname{Corollary #1}}

\def\myconjecture#1{\csname proclaim\endcsname{Conjecture #1}}
\def\myconjecturewithtitle#1#2{\csname proclaim\endcsname{Conjecture #1#2}}

\def\myproblem#1{\csname proclaim\endcsname{Problem #1}}
\def\myproblemwithtitle#1#2{\csname proclaim\endcsname{Problem #1#2}}
\def\mytheproblem#1{\blue{#1}\immediate\special{ps:
     ShrHPSdict begin /ShrBORDERthickness 0 def}}

\def\myanchortext#1#2{#2}
\def\mytheanchortext#1#2{\blue{#2}\immediate\special{ps:
     ShrHPSdict begin /ShrBORDERthickness 0 def}}
\font\eightcyr=wncyr8
\pagewidth{360pt}
\pageheight{606pt}
\topmatter
\title
On two algebraic parametrizations for rational solutions 
of the cuboid equations.
\endtitle
\rightheadtext{On two algebraic parametrizations \dots}
\author
John Ramsden, Ruslan Sharipov
\endauthor
\address CSR Plc, Cambridge Business Park, Cambridge, CB4 0WZ, UK
\endaddress
\email\myhref{mailto:jhnrmsdn\@yahoo.co.uk}{jhnrmsdn\@yahoo.co.uk}
\endemail
\address Bashkir State University, 32 Zaki Validi street, 450074 Ufa, Russia
\endaddress
\email\myhref{mailto:r-sharipov\@mail.ru}{r-sharipov\@mail.ru}
\endemail
\abstract
     A rational perfect cuboid is a rectangular parallelepiped whose edges and face 
diagonals are given by rational numbers and whose space diagonal is equal to unity. 
Its existence is equivalent to the existence of a perfect cuboid with all integer 
edges and diagonals. Finding such a cuboid or proving its non-existence is an old 
unsolved problem. Recently, based on a symmetry approach, the equations of a
perfect cuboid were transformed to factor equations. The factor equations turned 
out to be solvable and, being solved, have led to a pair of inverse problems. Our
efforts in the present paper are toward solving these inverse problems. Algebraic
parametrizations for their solutions using algebraic functions of two rational 
arguments are found. 
\endabstract
\endtopmatter
\TagsOnRight
\document

\head
1. Introduction.
\endhead
     Referring the reader to \myciterange{1}{1}{--}{44} for the history of cuboid
studies, we proceed to the following two cubic equations:
$$
\align
&\hskip -2em
x^3-E_{10}\,x^2+E_{20}\,x-E_{30}=0,
\mytag{1.1}\\
&\hskip -2em
d^{\kern 1pt 3}-E_{01}\,d^{\kern 1pt 2}+E_{02}\,d-E_{03}=0.
\mytag{1.2}
\endalign
$$
The equations \mythetag{1.1} and \mythetag{1.2} were derived as a result of the 
series of papers \myciterange{45}{45}{--}{50} along with three auxiliary equations 
$$
\hskip -2em
\aligned
&x_1\,x_2\,d_3+x_2\,x_3\,d_1+x_3\,x_1\,d_2=E_{21},\\
&x_1\,d_2+d_1\,x_2+x_2\,d_3+d_2\,x_3+x_3\,d_1+d_3\,x_1=E_{11},\\
&x_1\,d_2\,d_3+x_2\,d_3\,d_1+x_3\,d_1\,d_2=E_{12}.
\endaligned
\mytag{1.3}
$$
The numbers $x_1$, $x_2$, $x_3$ and $d_1$, $d_2$, $d_3$ in \mythetag{1.3} are 
edges and face diagonals of a rational perfect cuboid. The first three of them
are roots of the first cubic equation \mythetag{1.2}, the others are roots of 
the second cubic equation \mythetag{1.2}.\par
      The numbers $E_{10}$, $E_{20}$, $E_{30}$, $E_{01}$, $E_{02}$, $E_{03}$,
$E_{21}$, $E_{11}$, $E_{12}$ in the equations \mythetag{1.1}, \mythetag{1.2} and 
\mythetag{1.3} are rational functions of two arbitrary rational parameters $b$ 
and $c$. Here is the formula for the number $E_{11}$ in \mythetag{1.3}: 
$$
\hskip -2em
E_{11}=-\frac{b\,(c^2+2-4\,c)}{b^2\,c^2+2\,b^2-3\,b^2\,c+c-b\,c^2\,+2\,b}.
\mytag{1.4}
$$
The formulas for $E_{10}$, $E_{01}$ are similar to the formula \mythetag{1.4}
for $E_{11}$:
$$
\gather
\hskip -2em
E_{10}=-\frac{b^2\,c^2+2\,b^2-3\,b^2\,c\,-c}{b^2\,c^2+2\,b^2-3\,b^2\,c
+c-b\,c^2+2\,b},
\mytag{1.5}\\
\vspace{1ex}
\hskip -2em
E_{01}=-\frac{b\,(c^2+2-2\,c)}{b^2\,c^2+2\,b^2-3\,b^2\,c+c-b\,c^2+2\,b}.
\mytag{1.6}
\endgather
$$
Below are the formulas for $E_{20}$, $E_{02}$, $E_{30}$, 
$E_{03}$, $E_{21}$, $E_{12}$ in \mythetag{1.1}, \mythetag{1.2}, 
and \mythetag{1.3}:
$$
\gather
\hskip -2em
\gathered
E_{20}=\frac{b}{2}\,(b\,c^2-2\,c-2\,b)\,(2\,b\,c^2-c^2-6\,b\,c+2
+4\,b)\,\times\\
\times\,(b\,c-1-b)^{-2}\,(b\,c-c-2\,b)^{-2},
\endgathered\qquad\quad
\mytag{1.7}\\
\vspace{1ex}
\hskip -2em
\gathered
E_{02}=\frac{1}{2}\,(28\,b^2\,c^2-16\,b^2\,c-2\,c^2-4\,b^2-b^2\,c^4
+4\,b^3\,c^4-12\,b^3\,c^3\,+\\
+\,4\,b\,c^3+24\,b^3\,c-8\,b\,c-2\,b^4\,c^4+12\,b^4\,c^3-26\,b^4\,c^2
-8\,b^2\,c^3\,+\\
+24\,b^4\,c-16\,b^3-8\,b^4)\,(b\,c-1-b)^{-2}\,(b\,c-c-2\,b)^{-2},
\endgathered\qquad\quad
\mytag{1.8}\\
\vspace{1ex}
\hskip -2em
\gathered
E_{30}=c\,b^2\,(1-c)\,(c-2)\,(b\,c^2-4\,b\,c+2+4\,b)
\,(2\,b\,c^2-c^2-4\,b\,c\,+\\
+\,2\,b)\,(b^2\,c^4-6\,b^2\,c^3+13\,b^2\,c^2-12\,b^2\,c+4\,b^2
+c^2)^{-1}\,\times\\
\times\,(b\,c-1-b)^{-2}\,(-c+b\,c-2\,b)^{-2},
\endgathered\qquad\quad
\mytag{1.9}\\
\vspace{1ex}
\hskip -2em
\gathered
E_{03}=\frac{b}{2}\,(b^2\,c^4-5\,b^2\,c^3+10\,b^2\,c^2-10\,b^2\,c+4\,b^2
+2\,b\,c+2\,c^2\,-\\
-\,b\,c^3)\,(2\,b^2\,c^4-12\,b^2\,c^3+26\,b^2\,c^2-24\,b^2\,c
+\,8\,b^2-c^4\,b+3\,b\,c^3\,-\\
-\,6\,b\,c+4\,b+c^3-2\,c^2+2\,c)\,(b^2\,c^4-6\,b^2\,c^3+13\,b^2\,c^2\,-\\
-12\,b^2\,c+4\,b^2+c^2)^{-1}\,(b\,c-1-b)^{-2}\,(-c+b\,c-2\,b)^{-2},
\endgathered\qquad\quad
\mytag{1.10}\\
\vspace{1ex}
\gathered
E_{21}=\frac{b}{2}\,(5\,c^6\,b-2\,c^6\,b^2+52\,c^5\,b^2-16\,c^5\,b
-2\,c^7\,b^2+2\,b^4\,c^8\,-\\
-\,26\,b^4\,c^7-426\,b^4\,c^5-61\,b^3\,c^6+100\,b^3\,c^5
+14\,c^7\,b^3-c^8\,b^3-20\,b\,c^2\,-\\
-\,8\,b^2\,c^2-16\,b^2\,c-128\,b^2\,c^4-200\,b^3\,c^3
+244\,b^3\,c^2+32\,b\,c^3\,+\\
+\,768\,b^4\,c^4-852\,b^4\,c^3+568\,b^4\,c^2+104\,b^2\,c^3-208\,b^4\,c
+8\,c^4\,+\\
+16\,b^3-112\,b^3\,c+142\,b^4\,c^6
+32\,b^4-2\,c^5)\,(b^2\,c^4-6\,b^2\,c^3+13\,b^2\,c^2\,-\\
-12\,b^2\,c-4\,c^3+4\,b^2+c^2)^{-1}\,(b\,c-1-b)^{-2}\,(b\,c-c-2\,b)^{-2},
\endgathered\qquad\quad
\mytag{1.11}\\
\vspace{2ex}
\hskip -2em
\gathered
E_{12}=(16\,b^6+32\,b^5-6\,c^5\,b^2+2\,c^5\,b-62\,b^5\,c^6
+62\,b^6\,c^6+16\,b^4\,-\\
-\,180\,b^6\,c^5-c^7\,b^3+18\,b^5\,c^7-12\,b^6\,c^7-2\,b^5\,c^8
+b^6\,c^8+248\,b^5\,c^2\,+\\
+\,248\,b^6\,c^2-96\,b^6\,c+321\,b^6\,c^4-180\,b^5\,c^3-144\,b^5\,c
-360\,b^6\,c^3\,+\\
+\,b^4\,c^8+8\,b^4\,c^6-6\,b^4\,c^7+18\,b^4\,c^5+7\,b^3\,c^6
+90\,b^5\,c^5-14\,b^3\,c^5\,+\\
+\,17\,b^2\,c^4+32\,b^4\,c^2+28\,b^3\,c^3-28\,b^3\,c^2-4\,b\,c^3+8\,b^3\,c
-57\,b^4\,c^4\,+\\
+\,36\,b^4\,c^3-12\,b^2\,c^3-48\,b^4\,c-c^4)\,(b^2\,c^4-6\,b^2\,c^3
+13\,b^2\,c^2\,-\\
-\,12\,b^2\,c+4\,b^2+c^2)^{-1}\,(b\,c-1-b)^{-2}\,(b\,c-c-2\,b)^{-2}.
\endgathered\qquad
\mytag{1.12}
\endgather	
$$
Once the quantities $E_{10}$, $E_{20}$, $E_{30}$, $E_{01}$, $E_{02}$, $E_{03}$,
$E_{21}$, $E_{11}$, $E_{12}$ are known through the formulas \mythetag{1.4}, 
\mythetag{1.5}, \mythetag{1.6}, \mythetag{1.7}, \mythetag{1.8}, \mythetag{1.9}, 
\mythetag{1.10}, \mythetag{1.11}, \mythetag{1.12}, the next step is to find $x_1$, 
$x_2$, $x_3$ and $d_1$, $d_2$, $d_3$ by solving the equations \mythetag{1.1}, 
\mythetag{1.2}, \mythetag{1.3}. This step is formulated in the following inverse
problems. 
\myproblem{1.1} Find all pairs of rational numbers $b$ and $c$ for which the
cubic equations \mythetag{1.1} and \mythetag{1.2} with the coefficients given
by the formulas \mythetag{1.5}, \mythetag{1.7},	\mythetag{1.9} and \mythetag{1.6}, 
\mythetag{1.8},	\mythetag{1.10} possess positive rational roots $x_1$, $x_2$, 
$x_3$, $d_1$, $d_2$, $d_3$ obeying the auxiliary polynomial equations 
\mythetag{1.3} whose right hand sides are given by the formulas \mythetag{1.11}, 
\mythetag{1.4}, \mythetag{1.12}. 
\endproclaim
\myproblem{1.2} Find at least one pair of rational numbers $b$ and $c$ for which 
the cubic equations \mythetag{1.1} and \mythetag{1.2} with the coefficients given
by the formulas \mythetag{1.5}, \mythetag{1.7},	\mythetag{1.9} and \mythetag{1.6}, 
\mythetag{1.8},	\mythetag{1.10} possess positive rational roots $x_1$, $x_2$, 
$x_3$, $d_1$, $d_2$, $d_3$ obeying the auxiliary polynomial equations 
\mythetag{1.3} whose right hand sides are given by the formulas \mythetag{1.11}, 
\mythetag{1.4}, \mythetag{1.12}. 
\endproclaim
     The term ``inverse'' here means that $E_{10}$, $E_{20}$, $E_{30}$, $E_{01}$, 
$E_{02}$, $E_{03}$, $E_{21}$, $E_{11}$, $E_{12}$ are produced from $x_1$, 
$x_2$, $x_3$, and $d_1$, $d_2$, $d_3$ as the values of elementary multisymmetric 
polynomials (see \myciterange{51}{51}{--}{71}), which is treated as a direct
transform, then recovering $x_1$, $x_2$, $x_3$, $d_1$, $d_2$, $d_3$ through 
$E_{10}$, $E_{20}$, $E_{30}$, $E_{01}$, $E_{02}$, $E_{03}$, $E_{21}$, $E_{11}$, 
$E_{12}$ is an inverse transform. Since $E_{10}$, $E_{20}$, $E_{30}$, $E_{01}$, 
$E_{02}$, $E_{03}$, $E_{21}$, $E_{11}$, $E_{12}$ are functions of $b$ and $c$,
expressing $x_1$, $x_2$, $x_3$, $d_1$, $d_2$, $d_3$ through them means expressing
$x_1$, $x_2$, $x_3$, $d_1$, $d_2$, $d_3$ through $b$ and $c$ so that the equations
\mythetag{1.1}, \mythetag{1.2}, \mythetag{1.3} are fulfilled. Below in the present 
paper we find such expressions for $x_1$, $x_2$, $x_3$ and $d_1$, $d_2$, $d_3$ in
two ways using two algebraic function $w_1(b,c)$ and $w_2(b,c)$.\par
\head
2. Cubics with three rational roots. 
\endhead
\mylemma{2.1} A reduced cubic equation $y^3+y^2+D=0$ has three rational roots 
if and only if there is a rational number $w$ satisfying the sextic equation 
$$
\hskip -2em
D\,(w^2+3)^3+4\,(w-1)^2\,(1+w)^2=0.
\mytag{2.1}
$$
In this case the roots of the cubic equation $y^3+y^2+D=0$ are given by the
formulas 
$$
\xalignat 3
&\hskip -2em
y_1=-\frac{2\,(w+1)}{w^2+3},   
&&y_2=\frac{2\,(w-1)}{w^2+3},
&&y_3=\frac{1-w^2}{w^2+3}.
\quad
\mytag{2.2}
\endxalignat
$$
\endproclaim
\demo{Proof} Sufficiency. Assume that $w$ is a root of the sextic equation
\mythetag{2.1}. Note that the denominator $w^2+3$, which is common for all of
the three fractions \mythetag{2.2}, cannot vanish for any rational number $w$. 
Therefore \mythetag{2.2} yields three rational numbers $y_1$, $y_2$, $y_3$.
The rest is to substitute them into the product $(y-y_1)\,(y-y_2)\,(y-y_3)$:
$$
\hskip -2em
(y-y_1)\,(y-y_2)\,(y-y_3)=y^3+y^2-\frac{4\,(w-1)^2\,(1+w)^2}{(w^2+3)^3}.
\mytag{2.3}
$$
Since $w^2+3\neq 0$, the equation \mythetag{2.1} can be resolved with respect
to $D$:
$$
\hskip -2em
D=-\frac{4\,(w-1)^2\,(1+w)^2}{(w^2+3)^3}.
\mytag{2.4}
$$
Comparing \mythetag{2.4} and \mythetag{2.3}, we find that the sufficiency is
proved.\par
     Necessity. Assume that the equation $y^3+y^2+D=0$ has three rational roots
$y_1$, $y_2$, $y_3$. Then it can be written as 
$$
(y-y_2)\,(y^2+A\,y+B)=y^3+(A-y_2)\,y^2+(B-y_2\,A)\,y-y_2\,B=0,
\quad
\mytag{2.5}
$$
where $y_1$ and $y_3$ are roots of the quadratic equation $y^2+A\,y+B=0$. 
Comparing \mythetag{2.5} with the initial cubic equation $y^3+y^2+D=0$, we 
find
$$
\xalignat 3
&\hskip -2em
A=y_2+1,
&&B=y_2\,(y_2+1),
&&D=-y_2^3-y_2^2.
\mytag{2.6}
\endxalignat
$$
Note that the quadratic equation $y^2+A\,y+B=0$ with rational coefficients 
has two rational roots if and only if its discriminant is a square of some 
rational number $z$: 
$$
A^2-4\,B=z^2.
\mytag{2.7}
$$
Applying \mythetag{2.6} to \mythetag{2.7}, we derive the equation relating
$y_2$ and $z$:
$$
\hskip -2em
-3\,y_2^2-2\,y_2+1=z^2.
\mytag{2.8}
$$
The equation \mythetag{2.8} is similar to the equation (4.5) in \mycite{72}. 
It is solved similarly with the use of the lemma 2.2 in \mycite{49}. Its 
general solution in rational numbers is 
$$
\xalignat 2
&\hskip -2em
y_2=\frac{2\,t}{(t+1)^2+3},
&&z=\frac{t^2-4}{(t+1)^2+3},
\mytag{2.9}
\endxalignat
$$
where $t$ is an arbitrary rational number. The roots $y_1$ and $y_3$ of the 
quadratic equation $y^2+A\,y+B=0$ are given by the standard formula:
$$
\hskip -2em
y_{1,3}=-\frac{A}{2}\pm\frac{\sqrt{A^2-4\,B\vphantom{\vrule height 10pt}\,}}{2}. 
\mytag{2.10}
$$
Applying \mythetag{2.6} and \mythetag{2.7} to \mythetag{2.10}, we derive
$$
\xalignat 2
&\hskip -2em
y_1=-\frac{y_2+1}{2}+\frac{z}{2},
&&y_3=-\frac{y_2+1}{2}-\frac{z}{2},
\mytag{2.11}
\endxalignat
$$
Then we apply \mythetag{2.9} to \mythetag{2.11}. As a result we obtain
$$
\xalignat 3
&\hskip -2em
y_1=\frac{-2\,(2+t)}{(t+1)^2+3},
&&y_2=\frac{2\,t}{(t+1)^2+3},
&&y_3=\frac{-t\,(2+t)}{(t+1)^2+3}.
\qquad
\mytag{2.12}
\endxalignat
$$
In order to derive the required formulas \mythetag{2.2} it is sufficient to 
substitute $t=w-1$ into \mythetag{2.12}. Thus, we have found that if the 
cubic equation $y^3+y^2+D=0$ has three rational roots, these roots are 
expressed through some rational number $w$ by means of the formulas
\mythetag{2.2}. The rest is to substitute any one of the roots \mythetag{2.2}
into the initial cubic equation $y^3+y^2+D=0$. Since $w^2+3\neq 0$, this yields 
the required equation \mythetag{2.1}. The lemma~\mythelemma{2.1} is proved. 
\qed\enddemo
     Now let's consider a general cubic equation with the coefficients $A_0$, $A_1$, 
$A_2$, $A_3$:
$$
\hskip -2em
A_3\,x^3+A_2\,x^2+A_1\,x+A_0=0.
\mytag{2.13}
$$
Assuming that $A_3\neq 0$ in \mythetag{2.13}, we substitute
$$
\hskip -2em
x=\tilde x-\frac{A_2}{3\,A_3} 
\mytag{2.14}
$$
into the cubic equation \mythetag{2.13}. As a result it is transformed to
$$
\hskip -2em
\tilde x^3+\biggl(\frac{A_1}{A_3}-\frac{A_2^2}{3\,A_3^2}\biggr)\,\tilde x
+\frac{A_0}{A_3}-\frac{A_1\,A_2}{3\,A_3^2}+\frac{2\,A_2^3}{27\,A_3^3}=0.
\mytag{2.15}
$$
The following notations are for the sake of convenience:
$$
\xalignat 2
&\hskip -2em
B_1=\frac{A_1}{A_3}-\frac{A_2^2}{3\,A_3^2},
&&B_9=\frac{A_0}{A_3}-\frac{A_1\,A_2}{3\,A_3^2}+\frac{2\,A_2^3}{27\,A_3^3}.
\mytag{2.16}
\endxalignat
$$
In terms of the notations \mythetag{2.16} the cubic equation \mythetag{2.15}
simplifies to
$$
\hskip -2em
\tilde x^3+B_1\,\tilde x+B_0=0.
\mytag{2.17}
$$
Assume that $B_0\neq 0$ and $B_1\neq 0$. Then $\tilde x=0$ is not a root of the
cubic equation \mythetag{2.17}. Therefore the following transformation is applicable
to it:
$$
\hskip -2em
\tilde x=\frac{B_0}{B_1\,y}.
\mytag{2.18}
$$
Applying \mythetag{2.18} to the equation \mythetag{2.17}, we transform it to
$$
\hskip -2em
y^3+y^2+D=0\text{, \ where \ }D=\frac{B_0^2}{B_1^3}\neq 0. 
\mytag{2.19}
$$
Now, substituting \mythetag{2.16} into \mythetag{2.19}, we derive the formula 
for the parameter $D$ in the case of a general cubic equation \mythetag{2.13}:
$$
D=-\frac{(9\,A_1\,A_2\,A_3-27\,A_0\,A_3^2-2\,A_2^3)^2}
{27\,(A_2^2-3\,A_1\,A_3)^3}.
\mytag{2.20}
$$\par
     The next step is to transform the formulas \mythetag{2.2} for the roots 
of \mythetag{2.19} backward to the roots of a general cubic equation \mythetag{2.13} 
using \mythetag{2.18} and \mythetag{2.14}:
$$
\xalignat 3
&x_1=\frac{B_0}{B_1\,y_1}-\frac{A_2}{3\,A_3},
&&x_2=\frac{B_0}{B_1\,y_2}-\frac{A_2}{3\,A_3},
&&x_3=\frac{B_0}{B_1\,y_3}-\frac{A_2}{3\,A_3}.
\qquad\quad
\mytag{2.21}
\endxalignat
$$
Substituting \mythetag{2.16} into \mythetag{2.21} we derive the following three
formulas: 
$$
\allowdisplaybreaks
\gather
\hskip -2em
\gathered
x_1=\frac{1}{18}\,
((2\,A_2^3-9\,A_1\,A_2\,A_3+27\,A_0\,A_3^2)\,w^2
+(18\,A_2\,A_1\,A_3-6\,A_2^3)\,w\,-\\
-\,9\,A_1\,A_2\,A_3+81\,A_0\,A_3^2
)\,A_3^{-1}\,(A_2^2-3\,A_1\,A_3)^{-1}\,(1+w)^{-1},
\endgathered\quad
\mytag{2.22}\\
\displaybreak
\hskip -2em
\gathered
x_2=\frac{1}{18}\,
((2\,A_2^3-9\,A_1\,A_2\,A_3+27\,A_0\,A_3^2)\,w^2
-(18\,A_2\,A_1\,A_3-6\,A_2^3)\,w\,-\\
-\,9\,A_1\,A_2\,A_3+81\,A_0\,A_3^2
)\,A_3^{-1}\,(A_2^2-3\,A_1\,A_3)^{-1}\,(1-w)^{-1},
\endgathered\quad
\mytag{2.23}\\
\vspace{2ex}
\gathered
x_3=\frac{1}{9}\,((A_2^3-27\,A_0\,A_3^2)\,w^2+36\,A_1\,A_2\,A_3
-81\,A_0\,A_3^2-9\,A_2^3)\,\times\\
\times\,A_3^{-1}\,(A_2^2-3\,A_1\,A_3)^{-1}\,(1-w)^{-1}\,(1+w)^{-1}.
\endgathered\qquad
\mytag{2.24}
\endgather
$$
Based on the formulas \mythetag{2.22}, \mythetag{2.23}, \mythetag{2.24}, we can
formulate the next lemma. 
\mylemma{2.2} Assume that the numbers $A_0$, $A_1$, $A_2$, $A_3$ obey the 
inequalities
$$
\xalignat 3
&\hskip -2em
A_3\neq 0, 
&&\frac{A_1}{A_3}-\frac{A_2^2}{3\,A_3^2}\neq 0,
&&\frac{A_0}{A_3}-\frac{A_1\,A_2}{3\,A_3^2}+\frac{2\,A_2^3}{27\,A_3^3}\neq 0.
\qquad
\endxalignat
$$
Then the general cubic polynomial \mythetag{2.13} with the rational coefficients 
$A_0$, $A_1$, $A_2$, $A_3$ has three rational roots if and only if there is a 
rational number $w$ satisfying the sextic equation \mythetag{2.1} where $D$ is
given by the formula \mythetag{2.20}. In this case the roots of the cubic equation 
\mythetag{2.13} are given by the formulas \mythetag{2.22}, \mythetag{2.23},
\mythetag{2.24}. 
\endproclaim 
\head
3. The first algebraic parametrization.
\endhead
     The first algebraic parametrization for solutions of the equations
\mythetag{1.1}, \mythetag{1.2}, \mythetag{1.3} is produced from the the 
first cubic equation \mythetag{1.1} with the use of the lemma
\mythelemma{2.2}. Applying this lemma, we get the sextic equation 
$$
\hskip -2em
D_1\,(w^2+3)^3+4\,(w-1)^2\,(1+w)^2=0
\mytag{3.1}
$$
of the form \mythetag{2.1}. Its parameter $D=D_1$ is calculated using the formula
\mythetag{2.20} and the formulas \mythetag{1.5}, \mythetag{1.7}, \mythetag{1.9}
for the coefficients of the equation \mythetag{1.1}:
$$
\gathered
D_1=-\frac{2}{27}\,(7812\,b^4\,c^4\,-216\,b^2\,c^4-52\,b^2\,c^3+1764\,b^3\,c^4
-1200\,b^4\,c^3\,-\\
-\,1848\,b^4\,c^2+720\,b^4\,c-36\,c^4\,b-1512\,b^3\,c^3-36\,c^8\,b^3
+288\,b^3\,c^2\,-\\
-\,108\,c^6\,b^2+380\,c^5\,b^2+378\,c^7\,b^3-231\,c^8\,b^4-300\,c^7\,b^4
+3906\,c^6\,b^4\,-\\
-13\,c^7\,b^2-8904\,c^5\,b^4-882\,c^6\,b^3+18\,c^6\,b-1319\,b^6\,c^8
+20952\,b^5\,c^3\,-\\
-\,11952\,b^5\,c^2+2592\,b^5\,c-48372\,b^6\,c^4+31620\,b^6\,c^3-10552\,b^6\,c^2\,+\\
+\,\,816\,b^6\,c+1494\,b^5\,c^8-5238\,b^5\,c^7-4\,c^5+7905\,b^6\,c^7
-24186\,b^6\,c^6\,+\\
+\,288\,b^6+43740\,b^6\,c^5+7686\,b^5\,c^6+576\,b^7+128\,b^8-15372\,b^5\,c^4\,-\\
-\,1080\,b^7\,c^8-3546\,b^7\,c^6+51\,c^9\,b^6+400\,b^8\,c^8-162\,c^9\,b^5
+8640\,b^7\,c^2\,-\\
-\,3456\,b^7\,c+2808\,b^7\,c^7-1560\,b^8\,c^7+3940\,b^8\,c^6+216\,c^9\,b^7
-960\,b^8\,c\,-\\
-\,6240\,b^8\,c^3+9\,c^{10}\,b^6+7880\,b^8\,c^4+4\,c^{10}\,b^8-6732\,b^8\,c^5 
+45\,c^9\,b^4\,+\\
+\,3200\,b^8\,c^2-11232\,b^7\,c^3+7092\,b^7\,c^4-18\,c^{10}\,b^7
-60\,c^9\,b^8)^2\,(2\,c^2\,+\\
+\,2\,b^4\,c^4-12\,b^4\,c^3+26\,b^4\,c^2-24\,b^4\,c+8\,b^4-6\,b^3\,c^4
+18\,b^3\,c^3\,-\\
-\,36\,b^3\,c+24\,b^3+3\,b^2\,c^4+8\,b^2\,c^3-36\,b^2\,c^2+16\,b^2\,c+12\,b^2
-6\,b\,c^3\,+\\
+\,12\,b\,c)^3\,(b^2\,c^4-6\,b^2\,c^{-3}+13\,b^2\,c^2-12\,b^2\,c+4\,b^2+c^2)^{-2}.
\endgathered
\mytag{3.2}
$$
The fraction $2/27$ in \mythetag{3.2} can be expressed as $2^3/(2^2\,3^3)$. So the
structure of \mythetag{3.2} as a ratio of some square and some cube is the same as
the structure of $D$ in \mythetag{2.20}. The formulas \mythetag{3.1} and \mythetag{3.2} 
define an algebraic function $w=w_1(b,c)$. This function is used below as $w$ without 
showing its arguments.\par
     Now we proceed to the formulas \mythetag{2.22}, \mythetag{2.23}, \mythetag{2.24}.
Using these formulas, we get explicit expressions for $x_1$, $x_2$, $x_3$ through 
$b$, $c$, and $w$, where $w=w_1(b,c)$:
$$
\xalignat 3
&\hskip -2em
x_1=x_1(b,c,w),
&&x_2=x_2(b,c,w),
&&x_3=x_3(b,c,w).
\quad
\mytag{3.3}
\endxalignat
$$
However, the expressions for the functions \mythetag{3.3} are very huge. They 
comprise more than 100 terms in each. Therefore these expressions are given in 
\mytheanchortext{App1}{Appendix 1} in a machine readable form.\par
     The functions \mythetag{3.3} are roots of the cubic equation \mythetag{1.1}.
This fact follows from the lemma~\mythelemma{2.2}. Moreover, this fact has been 
tested computationally with the use of the explicit formulas for them in
\mytheanchortext{App1}{Appendix 1}.\par
     Apart from being roots of the cubic equation \mythetag{1.1}, the functions
\mythetag{3.3} obey the following cuboid equation saying that its space diagonal 
is equal to unity:
$$
\hskip -2em
x_1^2+x_2^2+x_3^2=1. 
\mytag{3.4}
$$
This fact follows from the theory of the cuboid factor equations in \mycite{46},
\mycite{47}. The equality \mythetag{3.4} has also been tested computationally. 
\par
     The nest step is to derive the formulas for the face diagonals $d_1$, $d_2$, 
$d_3$ of a cuboid. For this purpose we use the following equations:
$$
\hskip -2em
\aligned
&x_1\,x_2\,d_3+x_2\,x_3\,d_1+x_3\,x_1\,d_2=E_{21},\\
&x_1\,d_2+d_1\,x_2+x_2\,d_3+d_2\,x_3+x_3\,d_1+d_3\,x_1=E_{11},\\
&d_1+d_2+d_3=E_{01}.
\endaligned
\mytag{3.5}
$$
The first two of them are taken from \mythetag{1.3}. The last equation \mythetag{3.5}
follows from \mythetag{1.2} since $d_1$, $d_2$, $d_3$ should be the roots of this cubic 
equation.\par
     It is easy to see that the equations \mythetag{3.5} are linear with respect to 
$d_1$, $d_2$, $d_3$, while $x_1$, $x_2$, $x_3$ are already known from \mythetag{3.3}. 
Resolving the system of linear equations \mythetag{3.5}, we get three functions 
similar to $x_1$, $x_2$, $x_3$ in \mythetag{3.3}:
$$
\xalignat 3
&\hskip -2em
d_1=x_1(b,c,w),
&&d_2=d_2(b,c,w),
&&d_3=d_3(b,c,w).
\quad
\mytag{3.6}
\endxalignat
$$
The explicit formulas for the functions \mythetag{3.6} are extremely huge. They 
comprise more than 800 terms in each. For this reason we do not provide these 
formulas.\par
     The functions \mythetag{3.6} should be the roots of the second cubic equation
\mythetag{1.2}. Moreover, they should satisfy the last equation \mythetag{1.3}
which is not used in \mythetag{3.5} for determining them. These two facts follow
from the theory of the cuboid factor equations in \mycite{46}, \mycite{47}. 
However, they cannot be verified even with the use of symbolic computations since 
the formulas for the functions are extremely huge. These facts have been tested 
numerically for a series of random pairs of rational numbers $b$ and $c$.\par
\head
4. The second algebraic parametrization.
\endhead
    The second algebraic parametrization for solutions of the equations
\mythetag{1.1}, \mythetag{1.2}, \mythetag{1.3} is similar to the first one.
\pagebreak However, in tis case we start with the second cubic equation 
\mythetag{1.2}. Applying the lemma \mythelemma{2.2} to it, we get the 
sextic equation 
$$
\hskip -2em
D_2\,(w^2+3)^3+4\,(w-1)^2\,(1+w)^2=0
\mytag{4.1}
$$
of the form \mythetag{2.1}. Its parameter $D=D_2$ is calculated using the formula
\mythetag{2.20} and the formulas \mythetag{1.6}, \mythetag{1.8}, \mythetag{1.10}
for the coefficients of the equation \mythetag{1.2}:
$$
\gathered
D_2=-\frac{2\,b^2}{27}\,(832\,b^2\,c^2-1440\,b^2\,c^4-840\,b^2\,c^3
+4788\,b^3\,c^4+396\,b\,c^3\,+\\
+\,720\,b^3\,c+808\,b^4\,c^4+3032\,b^4\,c^3-2576\,b^4\,c^2
-96\,b^4\,c+448\,b^4\,-\\
-\,504\,c^4\,b-4176\,b^3\,c^3-9\,c^8\,b^3+72\,b^3\,c^2
-720\,c^6\,b^2+2288\,c^5\,b^2\,+\\
+\,1044\,c^7\,b^3-322\,c^8\,b^4+758\,c^7\,b^4+404\,c^6\,b^4
-210\,c^7\,b^2-2464\,c^5\,b^4\,-\\
-\,2394\,c^6\,b^3+72\,c^4+252\,c^6\,b+3168\,b^6\,c^8+441\,c^9\,b^5
-7056\,b^5\,c\,+\\
+\,57960\,b^6\,c^4-47232\,b^6\,c^3+25344\,b^6\,c^2
-8064\,b^6\,c-1809\,b^5\,c^8\,+\\
+\,14472\,b^5\,c^2+3951\,b^5\,c^7-72\,c^5+36\,c^6-11808\,b^6\,c^7
+1440\,b^5\,+\\
+\,28980\,b^6\,c^6-49032\,b^6\,c^5-4410\,b^5\,c^6
+8820\,b^5\,c^4-15804\,b^5\,c^3\,+\\
+\,1152\,b^6-504\,c^9\,b^6-45\,c^9\,b^3-6\,c^9\,b^4+104\,c^8\,b^2
+36\,c^{10}\,b^6\,+\\
+\,14\,c^{10}\,b^4-45\,c^{10}\,b^5-99\,c^7\,b)^2/(6\,b^4\,c^4
-36\,b^4\,c^3+78\,b^4\,c^2
-72\,b^4\,c\,+\\
+\,24\,b^4-12\,b^3\,c^4+36\,b^3\,c^3-72\,b^3\,c+48\,b^3+5\,b^2\,c^4
+16\,b^2\,c^3\,-\\
-\,68\,b^2\,c^2+32\,b^2\,c+20\,b^2-12\,b\,c^3+24\,b\,c
+6\,c^2)^{-3}\,(b^2\,c^4-6\,b^2\,c^3\,+\\
+\,13\,b^2\,c^2-12\,b^2\,c+4\,b^2+c^2)^{-2}.
\endgathered\quad
\mytag{4.2}
$$
The formulas \mythetag{4.1} and \mythetag{4.2} define an algebraic function $w=w_2(b,c)$. 
This function is used below as $w$ without showing its arguments.\par
     Now we proceed to the formulas \mythetag{2.22}, \mythetag{2.23}, \mythetag{2.24}.
Using these formulas, we get explicit expressions for $d_1$, $d_2$, $d_3$ through 
$b$, $c$, and $w$, where $w=w_2(b,c)$:
$$
\xalignat 3
&\hskip -2em
d_1=d_1(b,c,w),
&&d_2=d_2(b,c,w),
&&d_3=d_3(b,c,w).
\quad
\mytag{4.3}
\endxalignat
$$
Again, the expressions for the functions \mythetag{4.3} are very huge. They 
comprise more than 100 terms in each. Therefore these expressions are given in 
\mytheanchortext{App2}{Appendix 2} in a machine readable form.\par
     The functions \mythetag{4.3} are roots of the cubic equation \mythetag{1.2}.
This fact follows from the lemma~\mythelemma{2.2}. Moreover, this fact has been 
tested computationally with the use of the explicit formulas for them in
\mytheanchortext{App2}{Appendix 2}.\par
     Apart from being roots of the cubic equation \mythetag{1.1}, the functions
\mythetag{4.3} obey the following equation that can be derived from the original
cuboid equations:
$$
\hskip -2em
d_1^{\kern 1pt 2}+d_2^{\kern 1pt 2}+d_3^{\kern 1pt 2}=2. 
\mytag{4.4}
$$
However, which is more important, the equation \mythetag{4.4} follows from
the cuboid factor equations (see (6.10) in \mycite{46} or (1.12) and (1.19) in
\mycite{47} and recall \mythetag{3.4}). Despite the theoretical background from
\mycite{46} and \mycite{47}, the equality \mythetag{4.4} has also been tested 
computationally using the explicit expressions for \mythetag{4.3}.\par
     The nest step is to derive the formulas for $x_1$, $x_2$, $x_3$ through 
the formulas for $d_1$, $d_2$, $d_3$. For this purpose we use the following 
equations:
$$
\hskip -2em
\aligned
&x_1+x_2+x_3=E_{10},\\
&x_1\,d_2+d_1\,x_2+x_2\,d_3+d_2\,x_3+x_3\,d_1+d_3\,x_1=E_{11},\\
&x_1\,d_2\,d_3+x_2\,d_3\,d_1+x_3\,d_1\,d_2=E_{12}.
\endaligned
\mytag{4.5}
$$
The last two of them are taken from \mythetag{1.3}. The first equation \mythetag{4.5}
follows from \mythetag{1.1} since $x_1$, $x_2$, $x_3$ should be the roots of this cubic 
equation.\par
     It is easy to see that the equations \mythetag{4.5} are linear with respect to 
$x_1$, $x_2$, $x_3$, while $d_1$, $d_2$, $d_3$ are already known from \mythetag{4.3}. 
Resolving the system of linear equations \mythetag{4.5}, we get three functions 
similar to $d_1$, $d_2$, $d_3$ in \mythetag{4.3}:
$$
\xalignat 3
&\hskip -2em
x_1=x_1(b,c,w),
&&x_2=x_2(b,c,w),
&&x_3=x_3(b,c,w).
\quad
\mytag{4.6}
\endxalignat
$$
The explicit formulas for the functions \mythetag{4.6} are extremely huge. They 
comprise more than 800 terms in each. For this reason we do not provide these 
formulas.\par
     The functions \mythetag{4.6} should be the roots of the first cubic equation
\mythetag{1.1}. Moreover, they should satisfy the first equation \mythetag{1.3}
which is not used in \mythetag{4.5} for \nolinebreak determining them. These two 
facts follow from the theory of the cuboid factor equations in \mycite{46}, \mycite{47}. 
However, they cannot be verified even with the use of symbolic computations since 
the formulas for the functions are extremely huge. These facts have been tested 
numerically for a series of random pairs of rational numbers $b$ and $c$.\par
\head
5. Concluding remarks. 
\endhead
     Thus, two sets explicit formulas for possible solutions of the inverse cuboid
problems~\mytheproblem{1.1} and \mytheproblem{1.2} are obtained. However, in order
to produce an actual solution one should
solve at least one of the two sextic equations \mythetag{3.1} or \mythetag{4.1} in
rational numbers. The equations \mythetag{3.1} or \mythetag{4.1} are similar to
the twelfth order equation derived in \mycite{40}. They are sextic with respect to 
$w$. However their total degrees with respect to $b$, $c$ and $w$ are $42$ and $40$
respectively.\par 
     The equations \mythetag{3.1} or \mythetag{4.1} produce two algebraic functions
$w=w_1(b,c)$ and $w=w_2(b,c)$ which are different since the parameters $D_1$ and
$D_2$ in \mythetag{3.2} and \mythetag{4.2} are different. These two functions probably
are related to each other as $w_1=p_1(w_2,b,c)$ and $w_2=p_2(w_1,b,c)$, where $p_1$ and 
$p_2$ are polynomials in $w_2$ and $w_1$ respectively. But the relation can be 
more complicated, i\.\,e\. $P(w_1,w_2,b,c)=0$, where $P$ is a single polynomial of 
four variables. Which of these two options is valid? This question is to be studied 
in a separate paper.\par 
     The formulas \mythetag{3.1} and \mythetag{4.1} as well as the formulas in
\mytheanchortext{App1}{Appendix 1} and \mytheanchortext{App2}{Appendix} 
\mytheanchortext{App2}{2} and those huge formulas which are not presented explicitly 
have denominators. Some of them correspond to singularities studied in \mycite{73}. 
Others are new. These new singularities are also to be studied in a separate paper.
\par

\head
\myanchortext{App1}{Appendix 1.}
\endhead
Here are the formulas for $x_1=x_1(b,c,w)$, $x_2=x_2(b,c,w)$, $x_3=x_3(b,c,w)$
from \mythetag{3.3}. They are written in a machine readable form convenient 
for to copy-paste into some symbolic computations package:
\baselineskip=12pt plus 0.1pt minus 0.1pt
\medskip\noindent
{\tt
x1=1/18*(7686*w\^{}2*b\^{}5*c\^{}6+2592*b\^{}5*w\^{}2*c\relax
-36*w\^{}2*c\^{}8*b\^{}3+2808*w\^{}2*b\^{}7\newline
*c\^{}7+128*b\^{}8*w\^{}2-720*b\^{}2*c\^{}4-36*b\^{}2*c\^{}3\relax
+6084*b\^{}3*c\^{}4+26748*b\^{}4*c\^{}4-4\newline
176*b\^{}4*c\^{}3-6264*b\^{}4*c\^{}2+2448*b\^{}4*c\relax
+4*w\^{}2*c\^{}10*b\^{}8-36*c\^{}4*b-5400*b\^{}3\newline
*c\^{}3-1560*w\^{}2*b\^{}8*c\^{}7-144*c\^{}8*b\^{}3+\relax
1152*b\^{}3*c\^{}2-360*c\^{}6*b\^{}2+1116*c\^{}5*\newline
b\^{}2+1350*c\^{}7*b\^{}3-783*c\^{}8*b\^{}4-1044*c\^{}7*b\^{}4\relax
+13374*c\^{}6*b\^{}4-9*c\^{}7*b\^{}2-30\newline
456*c\^{}5*b\^{}4-3042*c\^{}6*b\^{}3+18*c\^{}6*b\relax
-4059*b\^{}6*c\^{}8+63720*b\^{}5*c\^{}3-36144*b\newline
\^{}5*c\^{}2+7776*b\^{}5*c-130716*b\^{}6*c\^{}4\relax
+88308*b\^{}6*c\^{}3-32472*b\^{}6*c\^{}2+4464*b\^{}\newline
6*c+4518*b\^{}5*c\^{}8-15930*b\^{}5*c\^{}7\relax
+9*w\^{}2*c\^{}10*b\^{}6+22077*b\^{}6*c\^{}7+288*b\^{}6-\newline
65358*b\^{}6*c\^{}6+117180*b\^{}6*c\^{}5+23454*b\^{}5*c\^{}6\relax
+31620*c\^{}3*b\^{}6*w\^{}2-6732*w\^{}\newline
2*b\^{}8*c\^{}5+1494*w\^{}2*b\^{}5*c\^{}8+576*b\^{}7\relax
-46908*b\^{}5*c\^{}4\relax
-1080*b\^{}7*c\^{}8+2808*b\newline
\^{}7*c\^{}7-3546*b\^{}7*c\^{}6+7092*b\^{}7*c\^{}4\relax
-11232*b\^{}7*c\^{}3+8640*b\^{}7*c\^{}2-3456*b\^{}7\newline
*c-486*c\^{}9*b\^{}5+279*c\^{}9*b\^{}6+216*c\^{}9*b\^{}7\relax
+153*c\^{}9*b\^{}4+9*c\^{}10*b\^{}6-18*c\^{}1\newline
0*b\^{}7-1319*w\^{}2*b\^{}6*c\^{}8-60*w\^{}2*c\^{}9*b\^{}8\relax
-1080*w\^{}2*b\^{}7*c\^{}8-162*w\^{}2*c\^{}9*b\newline
\^{}5+43740*w\^{}2*b\^{}6*c\^{}5+400*w\^{}2*b\^{}8*c\^{}8\relax
+51*w\^{}2*c\^{}9*b\^{}6-11232*w\^{}2*b\^{}7*c\^{}\newline
3-6240*w\^{}2*b\^{}8*c\^{}3-216*b\^{}2*w\^{}2*c\^{}4\relax
-11952*c\^{}2*b\^{}5*w\^{}2-3546*w\^{}2*b\^{}7*c\^{}\newline
6+3200*w\^{}2*b\^{}8*c\^{}2+3940*w\^{}2*b\^{}8*c\^{}6\relax
-300*w\^{}2*c\^{}7*b\^{}4+7905*w\^{}2*b\^{}6*c\^{}7\newline
+3906*w\^{}2*c\^{}6*b\^{}4+380*w\^{}2*c\^{}5*b\^{}2\relax
-1848*b\^{}4*w\^{}2*c\^{}2-5238*w\^{}2*b\^{}5*c\^{}7+\newline
7880*w\^{}2*b\^{}8*c\^{}4+816*c*b\^{}6*w\^{}2\relax
+576*b\^{}7*w\^{}2-231*w\^{}2*c\^{}8*b\^{}4-24186*w\^{}2\newline
*b\^{}6*c\^{}6+7092*w\^{}2*b\^{}7*c\^{}4\relax
+288*b\^{}6*w\^{}2-1152*b\^{}7*w-52*b\^{}2*w\^{}2*c\^{}3-4837\newline
2*w\^{}2*b\^{}6*c\^{}4+378*w\^{}2*c\^{}7*b\^{}3\relax
+216*w\^{}2*c\^{}9*b\^{}7+45*w\^{}2*c\^{}9*b\^{}4+1764*c\^{}\newline
4*b\^{}3*w\^{}2-720*b\^{}4*w*c\^{}2-72*b\^{}2*w*c\^{}4\relax
-864*b\^{}3*w*c\^{}3+720*b\^{}4*w\^{}2*c+288\newline
*b\^{}3*w\^{}2*c\^{}2+288*b\^{}4*w*c+288*b\^{}3*w*c\^{}2\relax
+120*b\^{}2*w*c\^{}3-4*w\^{}2*c\^{}5+72*b*\newline
w*c\^{}4-36*b*w\^{}2*c\^{}4-10552*c\^{}2*b\^{}6*w\^{}2\relax
-1512*b\^{}3*w\^{}2*c\^{}3+18*w\^{}2*c\^{}6*b-1\newline
8*c\^{}10*b\^{}7*w\^{}2-960*c*b\^{}8*w\^{}2\relax
-108*c\^{}6*b\^{}2*w\^{}2-3456*c*b\^{}7*w\^{}2-13*c\^{}7*b\newline
\^{}2*w\^{}2+8640*c\^{}2*b\^{}7*w\^{}2\relax
-882*c\^{}6*b\^{}3*w\^{}2+20952*c\^{}3*b\^{}5*w\^{}2-15372*c\^{}4*\newline
b\^{}5*w\^{}2+7812*c\^{}4*b\^{}4*w\^{}2\relax
-8904*c\^{}5*b\^{}4*w\^{}2-1200*c\^{}3*b\^{}4*w\^{}2-576*b\^{}6*w\newline
+3312*b\^{}4*w*c\^{}4-792*b\^{}5*w*c\^{}4\relax
+864*b\^{}5*w*c\^{}3-288*b\^{}5*w*c\^{}2-576*b\^{}4*w*\newline
c\^{}3-90*b\^{}4*w*c\^{}8-396*b\^{}3*w*c\^{}6\relax
-144*b\^{}4*w*c\^{}7-3744*b\^{}4*w*c\^{}5-24*b\^{}2*w\newline
*c\^{}5+36*b\^{}5*w*c\^{}8+396*b\^{}5*w*c\^{}6\relax
+1656*b\^{}4*w*c\^{}6-36*b\^{}2*w*c\^{}6+2160*b\^{}7\newline
|*w*c\^{}8-6552*b\^{}6*w*c\^{}3-816*b\^{}6*w*c\^{}2\relax
+2016*b\^{}6*w*c-14184*b\^{}7*w*c\^{}4+22\newline
464*b\^{}7*w*c\^{}3-17280*b\^{}7*w*c\^{}2-102*b\^{}6*w*c\^{}8\relax
-1638*b\^{}6*w*c\^{}7-14040*b\^{}6\newline
*w*c\^{}5-36*b*w*c\^{}6+7092*b\^{}7*w*c\^{}6+7200*b\^{}6*w*c\^{}6\relax
+14400*b\^{}6*w*c\^{}4-1200\newline
*b\^{}8*w*c\^{}8+4680*b\^{}8*w*c\^{}7-11820*b\^{}8*w*c\^{}6\relax
+20196*b\^{}8*w*c\^{}5+792*b\^{}3*w*\newline
c\^{}4-23640*b\^{}8*w*c\^{}4+18720*b\^{}8*w*c\^{}3\relax
-9600*b\^{}8*w*c\^{}2+126*b\^{}6*w*c\^{}9-432\newline
*b\^{}7*w*c\^{}9+180*b\^{}8*w*c\^{}9+30*b\^{}2*w*c\^{}7\relax
+18*b\^{}4*w*c\^{}9-36*b\^{}3*w*c\^{}8+36*b\newline
\^{}7*w*c\^{}10-12*b\^{}8*w*c\^{}10-18*b\^{}6*w*c\^{}10\relax
-384*b\^{}8*w+2880*w*b\^{}8*c+12*w*c\^{}\newline
5-5616*b\^{}7*w*c\^{}7+6912*b\^{}7*w*c-216*b\^{}5*w*c\^{}7\relax
+216*b\^{}3*w*c\^{}7)/((1+w)(20\newline
*b\^{}2*c\^{}2+16*b*c\^{}2+10*b\^{}2*c\^{}4+36*b\^{}2*c\relax
-66*b\^{}2*c\^{}3+78*b\^{}3*c\^{}4+80*b\^{}3*c\newline
-78*b\^{}4*c\^{}4+238*b\^{}4*c\^{}3-156*b\^{}4*c\^{}2\relax
-68*b\^{}4*c+24*b\^{}3+72*b\^{}4+2*c\^{}3-8*c\newline
\^{}4*b-156*b\^{}3*c\^{}2+9*c\^{}5*b\^{}2-20*c\^{}5*b\^{}3\relax
+9*c\^{}6*b\^{}4-17*c\^{}5*b\^{}4-3*c\^{}6*b\^{}3\newline
+176*b\^{}5*c\^{}2-192*b\^{}5*c+66*b\^{}6*c\^{}4\relax
-126*b\^{}6*c\^{}3+132*b\^{}6*c\^{}2-72*b\^{}6*c+4\newline
8*b\^{}5*c\^{}5+64*b\^{}5+16*b\^{}6+2*b\^{}6*c\^{}6\relax
-18*b\^{}6*c\^{}5-8*b\^{}5*c\^{}6-88*b\^{}5*c\^{}4)(b\newline
\^{}2*c\^{}4-6*b\^{}2*c\^{}3+13*b\^{}2*c\^{}2-12*b\^{}2*c\relax
+4*b\^{}2+c\^{}2));}
\medskip\noindent
{\tt
x2=-1/18*(7686*w\^{}2*b\^{}5*c\^{}6+2592*b\^{}5*w\^{}2*c\relax
-36*w\^{}2*c\^{}8*b\^{}3+2808*w\^{}2*b\^{}\newline
7*c\^{}7+128*b\^{}8*w\^{}2-720*b\^{}2*c\^{}4\relax
-36*b\^{}2*c\^{}3+6084*b\^{}3*c\^{}4+26748*b\^{}4*c\^{}4-\newline
4176*b\^{}4*c\^{}3-6264*b\^{}4*c\^{}2+2448*b\^{}4*c\relax
+4*w\^{}2*c\^{}10*b\^{}8-36*c\^{}4*b-5400*b\^{}\newline
3*c\^{}3-1560*w\^{}2*b\^{}8*c\^{}7-144*c\^{}8*b\^{}3\relax
+1152*b\^{}3*c\^{}2-360*c\^{}6*b\^{}2+1116*c\^{}5\newline
*b\^{}2+1350*c\^{}7*b\^{}3-783*c\^{}8*b\^{}4\relax
-1044*c\^{}7*b\^{}4+13374*c\^{}6*b\^{}4-9*c\^{}7*b\^{}2-3\newline
0456*c\^{}5*b\^{}4-3042*c\^{}6*b\^{}3+18*c\^{}6*b\relax
-4059*b\^{}6*c\^{}8+63720*b\^{}5*c\^{}3-36144*\newline
b\^{}5*c\^{}2+7776*b\^{}5*c-130716*b\^{}6*c\^{}4\relax
+88308*b\^{}6*c\^{}3-32472*b\^{}6*c\^{}2+4464*b\newline
\^{}6*c+4518*b\^{}5*c\^{}8-15930*b\^{}5*c\^{}7\relax
+9*w\^{}2*c\^{}10*b\^{}6+22077*b\^{}6*c\^{}7+288*b\^{}6\newline
-65358*b\^{}6*c\^{}6+117180*b\^{}6*c\^{}5+\relax
23454*b\^{}5*c\^{}6+31620*c\^{}3*b\^{}6*w\^{}2-6732*w\newline
\^{}2*b\^{}8*c\^{}5+1494*w\^{}2*b\^{}5*c\^{}8+\relax
576*b\^{}7-46908*b\^{}5*c\^{}4-1080*b\^{}7*c\^{}8+2808*\newline
b\^{}7*c\^{}7-3546*b\^{}7*c\^{}6+7092*b\^{}7*c\^{}4-\relax
11232*b\^{}7*c\^{}3+8640*b\^{}7*c\^{}2-3456*b\^{}\newline
7*c-486*c\^{}9*b\^{}5+279*c\^{}9*b\^{}6\relax
+216*c\^{}9*b\^{}7+153*c\^{}9*b\^{}4+9*c\^{}10*b\^{}6-18*c\^{}\newline
10*b\^{}7-1319*w\^{}2*b\^{}6*c\^{}8-60*w\^{}2*c\^{}9*b\^{}8\relax
-1080*w\^{}2*b\^{}7*c\^{}8-162*w\^{}2*c\^{}9*\newline
b\^{}5+43740*w\^{}2*b\^{}6*c\^{}5+400*w\^{}2*b\^{}8*c\^{}8\relax
+51*w\^{}2*c\^{}9*b\^{}6-11232*w\^{}2*b\^{}7*c\newline
\^{}3-6240*w\^{}2*b\^{}8*c\^{}3-216*b\^{}2*w\^{}2*c\^{}4\relax
-11952*c\^{}2*b\^{}5*w\^{}2-3546*w\^{}2*b\^{}7*c\newline
\^{}6+3200*w\^{}2*b\^{}8*c\^{}2+3940*w\^{}2*b\^{}8*c\^{}6\relax
-300*w\^{}2*c\^{}7*b\^{}4+7905*w\^{}2*b\^{}6*c\^{}\newline
7+3906*w\^{}2*c\^{}6*b\^{}4+380*w\^{}2*c\^{}5*b\^{}2\relax
-1848*b\^{}4*w\^{}2*c\^{}2-5238*w\^{}2*b\^{}5*c\^{}7\newline
+7880*w\^{}2*b\^{}8*c\^{}4+816*c*b\^{}6*w\^{}2\relax
+576*b\^{}7*w\^{}2-231*w\^{}2*c\^{}8*b\^{}4-24186*w\^{}\newline
2*b\^{}6*c\^{}6+7092*w\^{}2*b\^{}7*c\^{}4\relax
+288*b\^{}6*w\^{}2+1152*b\^{}7*w-52*b\^{}2*w\^{}2*c\^{}3-483\newline
72*w\^{}2*b\^{}6*c\^{}4+378*w\^{}2*c\^{}7*b\^{}3\relax
+216*w\^{}2*c\^{}9*b\^{}7+45*w\^{}2*c\^{}9*b\^{}4+1764*c\newline
\^{}4*b\^{}3*w\^{}2+720*b\^{}4*w*c\^{}2\relax
+72*b\^{}2*w*c\^{}4+864*b\^{}3*w*c\^{}3+720*b\^{}4*w\^{}2*c+28\newline
8*b\^{}3*w\^{}2*c\^{}2-288*b\^{}4*w*c\relax
-288*b\^{}3*w*c\^{}2-120*b\^{}2*w*c\^{}3-4*w\^{}2*c\^{}5-72*b\newline
*w*c\^{}4-36*b*w\^{}2*c\^{}4-10552*c\^{}2*b\^{}6*w\^{}2\relax
-1512*b\^{}3*w\^{}2*c\^{}3+18*w\^{}2*c\^{}6*b-\newline
18*c\^{}10*b\^{}7*w\^{}2-960*c*b\^{}8*w\^{}2\relax
-108*c\^{}6*b\^{}2*w\^{}2-3456*c*b\^{}7*w\^{}2-13*c\^{}7*\newline
b\^{}2*w\^{}2+8640*c\^{}2*b\^{}7*w\^{}2\relax
-882*c\^{}6*b\^{}3*w\^{}2+20952*c\^{}3*b\^{}5*w\^{}2-15372*c\^{}4\newline
*b\^{}5*w\^{}2+7812*c\^{}4*b\^{}4*w\^{}2\relax
-8904*c\^{}5*b\^{}4*w\^{}2-1200*c\^{}3*b\^{}4*w\^{}2+576*b\^{}6*\newline
w-3312*b\^{}4*w*c\^{}4+792*b\^{}5*w*c\^{}4\relax
-864*b\^{}5*w*c\^{}3+288*b\^{}5*w*c\^{}2+576*b\^{}4*w\newline
*c\^{}3+90*b\^{}4*w*c\^{}8+396*b\^{}3*w*c\^{}6\relax
+144*b\^{}4*w*c\^{}7+3744*b\^{}4*w*c\^{}5+24*b\^{}2*\newline
w*c\^{}5-36*b\^{}5*w*c\^{}8-396*b\^{}5*w*c\^{}6\relax
-1656*b\^{}4*w*c\^{}6+36*b\^{}2*w*c\^{}6-2160*b\^{}\newline
7*w*c\^{}8+6552*b\^{}6*w*c\^{}3\relax
+816*b\^{}6*w*c\^{}2-2016*b\^{}6*w*c+14184*b\^{}7*w*c\^{}4-22\newline
464*b\^{}7*w*c\^{}3+17280*b\^{}7*w*c\^{}2\relax
+102*b\^{}6*w*c\^{}8+1638*b\^{}6*w*c\^{}7+14040*b\^{}6\newline
*w*c\^{}5+36*b*w*c\^{}6-7092*b\^{}7*w*c\^{}6\relax
-7200*b\^{}6*w*c\^{}6-14400*b\^{}6*w*c\^{}4+1200\newline
*b\^{}8*w*c\^{}8-4680*b\^{}8*w*c\^{}7\relax
+11820*b\^{}8*w*c\^{}6-20196*b\^{}8*w*c\^{}5-792*b\^{}3*w*\newline
c\^{}4+23640*b\^{}8*w*c\^{}4-18720*b\^{}8*w*c\^{}3\relax
+9600*b\^{}8*w*c\^{}2-126*b\^{}6*w*c\^{}9+432\newline
*b\^{}7*w*c\^{}9-180*b\^{}8*w*c\^{}9-30*b\^{}2*w*c\^{}7\relax
-18*b\^{}4*w*c\^{}9+36*b\^{}3*w*c\^{}8-36*b\newline
\^{}7*w*c\^{}10+12*b\^{}8*w*c\^{}10+18*b\^{}6*w*c\^{}10\relax
+384*b\^{}8*w-2880*w*b\^{}8*c-12*w*c\^{}\newline
5+5616*b\^{}7*w*c\^{}7-6912*b\^{}7*w*c+216*b\^{}5*w*c\^{}7\relax
-216*b\^{}3*w*c\^{}7)/((20*b\^{}2*\newline
c\^{}2+16*b*c\^{}2+10*b\^{}2*c\^{}4+36*b\^{}2*c\relax
-66*b\^{}2*c\^{}3+78*b\^{}3*c\^{}4+80*b\^{}3*c-78*b\newline
\^{}4*c\^{}4+238*b\^{}4*c\^{}3-156*b\^{}4*c\^{}2\relax
-68*b\^{}4*c+24*b\^{}3+72*b\^{}4+2*c\^{}3-8*c\^{}4*b-\newline
156*b\^{}3*c\^{}2+9*c\^{}5*b\^{}2-20*c\^{}5*b\^{}3\relax
+9*c\^{}6*b\^{}4-17*c\^{}5*b\^{}4-3*c\^{}6*b\^{}3+176*\newline
b\^{}5*c\^{}2-192*b\^{}5*c+66*b\^{}6*c\^{}4-126*b\^{}6*c\^{}3\relax
+132*b\^{}6*c\^{}2-72*b\^{}6*c+48*b\^{}5\newline
*c\^{}5+64*b\^{}5+16*b\^{}6+2*b\^{}6*c\^{}6-18*b\^{}6*c\^{}5\relax
-8*b\^{}5*c\^{}6-88*b\^{}5*c\^{}4)(w-1)(b\newline
\^{}2*c\^{}4-6*b\^{}2*c\^{}3+13*b\^{}2*c\^{}2-12*b\^{}2*c\relax
+4*b\^{}2+c\^{}2));}
\medskip\noindent
{\tt 
x3=-2/9*(-3942*w\^{}2*b\^{}5*c\^{}6-1296*b\^{}5*w\^{}2*c\relax
+27*w\^{}2*c\^{}8*b\^{}3+32*b\^{}8*w\^{}2+\newline
306*b\^{}2*c\^{}4+108*b\^{}2*c\^{}3-2448*b\^{}3*c\^{}4-10890*b\^{}4*c\^{}4\relax
+1656*b\^{}4*c\^{}3+2592\newline
*b\^{}4*c\^{}2-1008*b\^{}4*c+w\^{}2*c\^{}10*b\^{}8+72*c\^{}4*b\relax
+2052*b\^{}3*c\^{}3-390*w\^{}2*b\^{}8*c\newline
\^{}7+45*c\^{}8*b\^{}3-360*b\^{}3*c\^{}2+153*c\^{}6*b\^{}2\relax
-576*c\^{}5*b\^{}2-513*c\^{}7*b\^{}3+324*c\^{}\newline
8*b\^{}4+414*c\^{}7*b\^{}4-5445*c\^{}6*b\^{}4+27*c\^{}7*b\^{}2\relax
+12420*c\^{}5*b\^{}4+1224*c\^{}6*b\^{}3\newline
-36*c\^{}6*b+1953*b\^{}6*c\^{}8-31212*b\^{}5*c\^{}3+17856*b\^{}5*c\^{}2\relax
-3888*b\^{}5*c+76158*\newline
b\^{}6*c\^{}4-49068*b\^{}6*c\^{}3+15624*b\^{}6*c\^{}2-720*b\^{}6*c\relax
-2232*b\^{}5*c\^{}8+7803*b\^{}5*\newline
c\^{}7+9*c\^{}5-12267*b\^{}6*c\^{}7-576*b\^{}6+38079*b\^{}6*c\^{}6\relax
-69120*b\^{}6*c\^{}5-11430*b\^{}\newline
5*c\^{}6-14172*c\^{}3*b\^{}6*w\^{}2-1683*w\^{}2*b\^{}8*c\^{}5\relax
-756*w\^{}2*b\^{}5*c\^{}8-1152*b\^{}7-28\newline
8*b\^{}8+22860*b\^{}5*c\^{}4+2160*b\^{}7*c\^{}8-5616*b\^{}7*c\^{}7\relax
+7092*b\^{}7*c\^{}6-14184*b\^{}7\newline
*c\^{}4+22464*b\^{}7*c\^{}3-17280*b\^{}7*c\^{}2+6912*b\^{}7*c\relax
-900*b\^{}8*c\^{}8+3510*b\^{}8*c\^{}7\newline
-8865*b\^{}8*c\^{}6+15147*b\^{}8*c\^{}5-17730*b\^{}8*c\^{}4\relax
+14040*b\^{}8*c\^{}3-7200*b\^{}8*c\^{}2\newline
+2160*b\^{}8*c+243*c\^{}9*b\^{}5-45*c\^{}9*b\^{}6-432*c\^{}9*b\^{}7\relax
-63*c\^{}9*b\^{}4-18*c\^{}10*b\^{}\newline
6+36*c\^{}10*b\^{}7-9*c\^{}10*b\^{}8+135*c\^{}9*b\^{}8\relax
+685*w\^{}2*b\^{}6*c\^{}8-15*w\^{}2*c\^{}9*b\^{}8+\newline
81*w\^{}2*c\^{}9*b\^{}5-18360*w\^{}2*b\^{}6*c\^{}5\relax
+100*w\^{}2*b\^{}8*c\^{}8-57*w\^{}2*c\^{}9*b\^{}6-1560\newline
*w\^{}2*b\^{}8*c\^{}3+126*b\^{}2*w\^{}2*c\^{}4\relax
+6048*c\^{}2*b\^{}5*w\^{}2+800*w\^{}2*b\^{}8*c\^{}2+985*w\^{}\newline
2*b\^{}8*c\^{}6+186*w\^{}2*c\^{}7*b\^{}4\relax
-3543*w\^{}2*b\^{}6*c\^{}7-2367*w\^{}2*c\^{}6*b\^{}4-184*w\^{}2*\newline
c\^{}5*b\^{}2+1104*b\^{}4*w\^{}2*c\^{}2+2673*w\^{}2*b\^{}5*c\^{}7\relax
+1970*w\^{}2*b\^{}8*c\^{}4-912*c*b\^{}6\newline
*w\^{}2+138*w\^{}2*c\^{}8*b\^{}4+10293*w\^{}2*b\^{}6*c\^{}6\relax
-4*b\^{}2*w\^{}2*c\^{}3+20586*w\^{}2*b\^{}6*c\newline
\^{}4-243*w\^{}2*c\^{}7*b\^{}3-27*w\^{}2*c\^{}9*b\^{}4\relax
-1080*c\^{}4*b\^{}3*w\^{}2-432*b\^{}4*w\^{}2*c-216\newline
*b\^{}3*w\^{}2*c\^{}2-w\^{}2*c\^{}5+5480*c\^{}2*b\^{}6*w\^{}2\relax
+972*b\^{}3*w\^{}2*c\^{}3-240*c*b\^{}8*w\^{}2+\newline
63*c\^{}6*b\^{}2*w\^{}2-c\^{}7*b\^{}2*w\^{}2\relax
+540*c\^{}6*b\^{}3*w\^{}2-10692*c\^{}3*b\^{}5*w\^{}2+7884*c\^{}\newline
4*b\^{}5*w\^{}2-4734*c\^{}4*b\^{}4*w\^{}2\relax
+5388*c\^{}5*b\^{}4*w\^{}2+744*c\^{}3*b\^{}4*w\^{}2)/((1+w)*\newline
(w-1)*(20*b\^{}2*c\^{}2+16*b*c\^{}2+10*b\^{}2*c\^{}4\relax
+36*b\^{}2*c-66*b\^{}2*c\^{}3+78*b\^{}3*c\^{}4\newline
+80*b\^{}3*c-78*b\^{}4*c\^{}4+238*b\^{}4*c\^{}3\relax
-156*b\^{}4*c\^{}2-68*b\^{}4*c+24*b\^{}3+72*b\^{}4+\newline
2*c\^{}3-8*c\^{}4*b-156*b\^{}3*c\^{}2+9*c\^{}5*b\^{}2\relax
-20*c\^{}5*b\^{}3+9*c\^{}6*b\^{}4-17*c\^{}5*b\^{}4-\newline
3*c\^{}6*b\^{}3+176*b\^{}5*c\^{}2-192*b\^{}5*c\relax
+66*b\^{}6*c\^{}4-126*b\^{}6*c\^{}3+132*b\^{}6*c\^{}2-7\newline
2*b\^{}6*c+48*b\^{}5*c\^{}5+64*b\^{}5+16*b\^{}6\relax
+2*b\^{}6*c\^{}6-18*b\^{}6*c\^{}5-8*b\^{}5*c\^{}6-88*b\newline
\^{}5*c\^{}4)*(b\^{}2*c\^{}4-6*b\^{}2*c\^{}3+13*b\^{}2*c\^{}2\relax
-12*b\^{}2*c+4*b\^{}2+c\^{}2));}\par

\head
\myanchortext{App2}{Appendix 2.}
\endhead
Here are the formulas for $d_1=d_1(b,c,w)$, $d_2=d_2(b,c,w)$, $d_3=d_3(b,c,w)$
from \mythetag{4.3}. They are written in a machine readable form convenient 
for to copy-paste into some symbolic computations package:
\vskip 0 pt plus 1 pt minus 1pt
\medskip\noindent
{\tt 
d1=-1/18*(-4410*w\^{}2*b\^{}5*c\^{}6+30*c\^{}10*b\^{}4*w\relax
+792*c\^{}9*b\^{}5*w-45*w\^{}2*c\^{}9*b\newline
\^{}3-7056*b\^{}5*w\^{}2*c-9*w\^{}2*c\^{}8*b\^{}3+3024*b\^{}2*c\^{}2\relax
-3456*b\^{}2*c\^{}4-3528*b\^{}2*c\newline
\^{}3+10620*b\^{}3*c\^{}4+1476*b*c\^{}3+3312*b\^{}3*c\relax
-26856*b\^{}4*c\^{}4+26136*b\^{}4*c\^{}3-1\newline
0368*b\^{}4*c\^{}2-2592*b\^{}4*c+2304*b\^{}4-1800*c\^{}4*b\relax
-6768*b\^{}3*c\^{}3+477*c\^{}8*b\^{}3\newline
-3816*b\^{}3*c\^{}2-1728*c\^{}6*b\^{}2+6264*c\^{}5*b\^{}2\relax
+1692*c\^{}7*b\^{}3-1296*c\^{}8*b\^{}4+65\newline
34*c\^{}7*b\^{}4-13428*c\^{}6*b\^{}4-882*c\^{}7*b\^{}2\relax
+19656*c\^{}5*b\^{}4-5310*c\^{}6*b\^{}3+900*\newline
c\^{}6*b+12672*b\^{}6*c\^{}8-83412*b\^{}5*c\^{}3+72792*b\^{}5*c\^{}2\relax
-33840*b\^{}5*c+231840*b\newline
\^{}6*c\^{}4-188928*b\^{}6*c\^{}3+101376*b\^{}6*c\^{}2\relax
-32256*b\^{}6*c-9099*b\^{}5*c\^{}8+20853*\newline
b\^{}5*c\^{}7+288*c\^{}4-288*c\^{}5+144*c\^{}6\relax
+36*w\^{}2*c\^{}10*b\^{}6-47232*b\^{}6*c\^{}7+6624*b\newline
\^{}5+4608*b\^{}6+115920*b\^{}6*c\^{}6-196128*b\^{}6*c\^{}5\relax
-24174*b\^{}5*c\^{}6-47232*c\^{}3*b\^{}\newline
6*w\^{}2-1809*w\^{}2*b\^{}5*c\^{}8+104*w\^{}2*c\^{}8*b\^{}2\relax
+48348*b\^{}5*c\^{}4+2115*c\^{}9*b\^{}5-20\newline
16*c\^{}9*b\^{}6-207*c\^{}9*b\^{}3-162*c\^{}9*b\^{}4\relax
+378*c\^{}8*b\^{}2+144*c\^{}10*b\^{}6+72*c\^{}10*\newline
b\^{}4+3168*w\^{}2*b\^{}6*c\^{}8+441*w\^{}2*c\^{}9*b\^{}5\relax
-49032*w\^{}2*b\^{}6*c\^{}5+72*w\^{}2*c\^{}4-72\newline
*c\^{}9*b\^{}3*w-504*w\^{}2*c\^{}9*b\^{}6-1440*b\^{}2*w\^{}2*c\^{}4\relax
+14472*c\^{}2*b\^{}5*w\^{}2+758*w\^{}\newline
2*c\^{}7*b\^{}4+36*w*c\^{}6-11808*w\^{}2*b\^{}6*c\^{}7\relax
+404*w\^{}2*c\^{}6*b\^{}4+2288*w\^{}2*c\^{}5*b\^{}\newline
2-2576*b\^{}4*w\^{}2*c\^{}2+3951*w\^{}2*b\^{}5*c\^{}7\relax
+72*w*c\^{}4-8064*c*b\^{}6*w\^{}2-322*w\^{}2*\newline
c\^{}8*b\^{}4+28980*w\^{}2*b\^{}6*c\^{}6+1152*b\^{}6*w\^{}2\relax
+36*w\^{}2*c\^{}6-840*b\^{}2*w\^{}2*c\^{}3+14\newline
*w\^{}2*c\^{}10*b\^{}4+57960*w\^{}2*b\^{}6*c\^{}4\relax
+1044*w\^{}2*c\^{}7*b\^{}3-6*w\^{}2*c\^{}9*b\^{}4+4788*\newline
c\^{}4*b\^{}3*w\^{}2-2640*b\^{}4*w*c\^{}2+864*b\^{}2*w*c\^{}4\relax
+5760*b\^{}3*w*c\^{}3-96*b\^{}4*w\^{}2*c\newline
-12672*b\^{}5*c*w+72*b\^{}3*w\^{}2*c\^{}2+1440*b\^{}5*w\^{}2\relax
-2304*b\^{}4*w*c-72*c\^{}10*b\^{}5*\newline
w+2304*b\^{}5*w-4032*b\^{}3*w*c\^{}2-1008*b\^{}2*w*c\^{}3\relax
-45*b\^{}5*w\^{}2*c\^{}10-72*w\^{}2*c\^{}\newline
5-288*b*w*c\^{}4-504*b*w\^{}2*c\^{}4+66*c\^{}8*b\^{}2*w\relax
-207*c\^{}10*b\^{}5+25344*c\^{}2*b\^{}6*\newline
w\^{}2-4176*b\^{}3*w\^{}2*c\^{}3+252*w\^{}2*c\^{}6*b\relax
-720*c\^{}6*b\^{}2*w\^{}2-210*c\^{}7*b\^{}2*w\^{}2-2\newline
394*c\^{}6*b\^{}3*w\^{}2-15804*c\^{}3*b\^{}5*w\^{}2\relax
+8820*c\^{}4*b\^{}5*w\^{}2+808*c\^{}4*b\^{}4*w\^{}2-2\newline
464*c\^{}5*b\^{}4*w\^{}2+3032*c\^{}3*b\^{}4*w\^{}2\relax
+1152*b\^{}6*w-29280*b\^{}4*w*c\^{}4+21888*b\^{}\newline
5*w*c\^{}4-36000*b\^{}5*w*c\^{}3+29376*b\^{}5*w*c\^{}2\relax
+17040*b\^{}4*w*c\^{}3-330*b\^{}4*w*c\^{}\newline
8+1872*b\^{}3*w*c\^{}6+4260*b\^{}4*w*c\^{}7\relax
+27048*b\^{}4*w*c\^{}5-600*b\^{}2*w*c\^{}5-3672*b\newline
\^{}5*w*c\^{}8-10944*b\^{}5*w*c\^{}6-14640*b\^{}4*w*c\^{}6\relax
+432*b\^{}2*w*c\^{}6-47232*b\^{}6*w*c\newline
\^{}3+25344*b\^{}6*w*c\^{}2-8064*b\^{}6*w*c\relax
+3168*b\^{}6*w*c\^{}8-11808*b\^{}6*w*c\^{}7-49032\newline
*b\^{}6*w*c\^{}5+144*b*w*c\^{}6+28980*b\^{}6*w*c\^{}6\relax
+57960*b\^{}6*w*c\^{}4-3744*b\^{}3*w*c\^{}\newline
4-504*b\^{}6*w*c\^{}9-252*b\^{}2*w*c\^{}7\relax
-144*b\^{}4*w*c\^{}9+504*b\^{}3*w*c\^{}8+36*b\^{}6*w*c\newline
\^{}10-72*w*c\^{}5+9000*b\^{}5*w*c\^{}7\relax
-1440*b\^{}3*w*c\^{}7+960*b\^{}4*w+448*b\^{}4*w\^{}2-99*\newline
b*w\^{}2*c\^{}7-72*b*w*c\^{}7-369*c\^{}7*b\relax
+720*b\^{}3*w\^{}2*c+1152*b\^{}3*c*w+528*b\^{}2*w*\newline
c\^{}2+832*b\^{}2*w\^{}2*c\^{}2+396*b*w\^{}2*c\^{}3\relax
+288*b*w*c\^{}3)*b/((1+w)(44*b\^{}2*c\^{}2+3\newline
6*b*c\^{}2+22*b\^{}2*c\^{}4+68*b\^{}2*c\relax
-134*b\^{}2*c\^{}3+150*b\^{}3*c\^{}4+160*b\^{}3*c-166*b\^{}\newline
4*c\^{}4+490*b\^{}4*c\^{}3-332*b\^{}4*c\^{}2\relax
-116*b\^{}4*c+40*b\^{}3+136*b\^{}4+6*c\^{}3-18*c\^{}4*\newline
b-300*b\^{}3*c\^{}2+17*c\^{}5*b\^{}2-40*c\^{}5*b\^{}3\relax
+17*c\^{}6*b\^{}4-29*c\^{}5*b\^{}4-5*c\^{}6*b\^{}3+\newline
396*b\^{}5*c\^{}2-432*b\^{}5*c+198*b\^{}6*c\^{}4\relax
-378*b\^{}6*c\^{}3+396*b\^{}6*c\^{}2-216*b\^{}6*c+\newline
108*b\^{}5*c\^{}5+144*b\^{}5+48*b\^{}6+6*b\^{}6*c\^{}6\relax
-54*b\^{}6*c\^{}5-18*b\^{}5*c\^{}6-198*b\^{}5*c\newline
\^{}4)(b\^{}2*c\^{}4-6*b\^{}2*c\^{}3+13*b\^{}2*c\^{}2\relax
-12*b\^{}2*c+4*b\^{}2+c\^{}2));}
\medskip\noindent
{\tt 
d2=1/18*(-4410*w\^{}2*b\^{}5*c\^{}6-30*c\^{}10*b\^{}4*w\relax
-792*c\^{}9*b\^{}5*w-45*w\^{}2*c\^{}9*b\^{}\newline
3-7056*b\^{}5*w\^{}2*c-9*w\^{}2*c\^{}8*b\^{}3\relax
+3024*b\^{}2*c\^{}2-3456*b\^{}2*c\^{}4-3528*b\^{}2*c\^{}\newline
3+10620*b\^{}3*c\^{}4+1476*b*c\^{}3+3312*b\^{}3*c\relax
-26856*b\^{}4*c\^{}4+26136*b\^{}4*c\^{}3-10\newline
368*b\^{}4*c\^{}2-2592*b\^{}4*c+2304*b\^{}4-1800*c\^{}4*b\relax
-6768*b\^{}3*c\^{}3+477*c\^{}8*b\^{}3-\newline
3816*b\^{}3*c\^{}2-1728*c\^{}6*b\^{}2+6264*c\^{}5*b\^{}2\relax
+1692*c\^{}7*b\^{}3-1296*c\^{}8*b\^{}4+653\newline
4*c\^{}7*b\^{}4-13428*c\^{}6*b\^{}4-882*c\^{}7*b\^{}2\relax
+19656*c\^{}5*b\^{}4-5310*c\^{}6*b\^{}3+900*c\newline
\^{}6*b+12672*b\^{}6*c\^{}8-83412*b\^{}5*c\^{}3\relax
+72792*b\^{}5*c\^{}2-33840*b\^{}5*c+231840*b\^{}\newline
6*c\^{}4-188928*b\^{}6*c\^{}3+101376*b\^{}6*c\^{}2\relax
-32256*b\^{}6*c-9099*b\^{}5*c\^{}8+20853*b\newline
\^{}5*c\^{}7+288*c\^{}4-288*c\^{}5+144*c\^{}6\relax
+36*w\^{}2*c\^{}10*b\^{}6-47232*b\^{}6*c\^{}7+6624*b\^{}\newline
5+4608*b\^{}6+115920*b\^{}6*c\^{}6-196128*b\^{}6*c\^{}5\relax
-24174*b\^{}5*c\^{}6-47232*c\^{}3*b\^{}6\newline
*w\^{}2-1809*w\^{}2*b\^{}5*c\^{}8+104*w\^{}2*c\^{}8*b\^{}2\relax
+48348*b\^{}5*c\^{}4+2115*c\^{}9*b\^{}5-201\newline
6*c\^{}9*b\^{}6-207*c\^{}9*b\^{}3-162*c\^{}9*b\^{}4\relax
+378*c\^{}8*b\^{}2+144*c\^{}10*b\^{}6+72*c\^{}10*b\newline
\^{}4+3168*w\^{}2*b\^{}6*c\^{}8+441*w\^{}2*c\^{}9*b\^{}5\relax
-49032*w\^{}2*b\^{}6*c\^{}5+72*w\^{}2*c\^{}4+72*\newline
c\^{}9*b\^{}3*w-504*w\^{}2*c\^{}9*b\^{}6\relax
-1440*b\^{}2*w\^{}2*c\^{}4+14472*c\^{}2*b\^{}5*w\^{}2+758*w\^{}2\newline
*c\^{}7*b\^{}4-36*w*c\^{}6-11808*w\^{}2*b\^{}6*c\^{}7\relax
+404*w\^{}2*c\^{}6*b\^{}4+2288*w\^{}2*c\^{}5*b\^{}2\newline
-2576*b\^{}4*w\^{}2*c\^{}2+3951*w\^{}2*b\^{}5*c\^{}7-72*w*c\^{}4\relax
-8064*c*b\^{}6*w\^{}2-322*w\^{}2*c\newline
\^{}8*b\^{}4+28980*w\^{}2*b\^{}6*c\^{}6\relax
+1152*b\^{}6*w\^{}2+36*w\^{}2*c\^{}6-840*b\^{}2*w\^{}2*c\^{}3+14*\newline
w\^{}2*c\^{}10*b\^{}4+57960*w\^{}2*b\^{}6*c\^{}4\relax
+1044*w\^{}2*c\^{}7*b\^{}3-6*w\^{}2*c\^{}9*b\^{}4+4788*c\newline
\^{}4*b\^{}3*w\^{}2+2640*b\^{}4*w*c\^{}2-864*b\^{}2*w*c\^{}4\relax
-5760*b\^{}3*w*c\^{}3-96*b\^{}4*w\^{}2*c+\newline
12672*b\^{}5*c*w+72*b\^{}3*w\^{}2*c\^{}2+1440*b\^{}5*w\^{}2\relax
+2304*b\^{}4*w*c+72*c\^{}10*b\^{}5*w\newline
-2304*b\^{}5*w+4032*b\^{}3*w*c\^{}2+1008*b\^{}2*w*c\^{}3\relax
-45*b\^{}5*w\^{}2*c\^{}10-72*w\^{}2*c\^{}5\newline
+288*b*w*c\^{}4-504*b*w\^{}2*c\^{}4-66*c\^{}8*b\^{}2*w\relax
-207*c\^{}10*b\^{}5+25344*c\^{}2*b\^{}6*w\newline
\^{}2-4176*b\^{}3*w\^{}2*c\^{}3+252*w\^{}2*c\^{}6*b\relax
-720*c\^{}6*b\^{}2*w\^{}2-210*c\^{}7*b\^{}2*w\^{}2-23\newline
94*c\^{}6*b\^{}3*w\^{}2-15804*c\^{}3*b\^{}5*w\^{}2\relax
+8820*c\^{}4*b\^{}5*w\^{}2+808*c\^{}4*b\^{}4*w\^{}2-24\newline
64*c\^{}5*b\^{}4*w\^{}2+3032*c\^{}3*b\^{}4*w\^{}2\relax
-1152*b\^{}6*w+29280*b\^{}4*w*c\^{}4-21888*b\^{}5\newline
*w*c\^{}4+36000*b\^{}5*w*c\^{}3-29376*b\^{}5*w*c\^{}2\relax
-17040*b\^{}4*w*c\^{}3+330*b\^{}4*w*c\^{}8\newline
-1872*b\^{}3*w*c\^{}6-4260*b\^{}4*w*c\^{}7\relax
-27048*b\^{}4*w*c\^{}5+600*b\^{}2*w*c\^{}5+3672*b\^{}\newline
5*w*c\^{}8+10944*b\^{}5*w*c\^{}6+14640*b\^{}4*w*c\^{}6\relax
-432*b\^{}2*w*c\^{}6+47232*b\^{}6*w*c\^{}\newline
3-25344*b\^{}6*w*c\^{}2+8064*b\^{}6*w*c-3168*b\^{}6*w*c\^{}8\relax
+11808*b\^{}6*w*c\^{}7+49032*\newline
b\^{}6*w*c\^{}5-144*b*w*c\^{}6-28980*b\^{}6*w*c\^{}6\relax
-57960*b\^{}6*w*c\^{}4+3744*b\^{}3*w*c\^{}4\newline
+504*b\^{}6*w*c\^{}9+252*b\^{}2*w*c\^{}7+144*b\^{}4*w*c\^{}9\relax
-504*b\^{}3*w*c\^{}8-36*b\^{}6*w*c\^{}\newline
10+72*w*c\^{}5-9000*b\^{}5*w*c\^{}7+1440*b\^{}3*w*c\^{}7\relax
-960*b\^{}4*w+448*b\^{}4*w\^{}2-99*b\newline
*w\^{}2*c\^{}7+72*b*w*c\^{}7-369*c\^{}7*b+720*b\^{}3*w\^{}2*c\relax
-1152*b\^{}3*c*w-528*b\^{}2*w*c\newline
\^{}2+832*b\^{}2*w\^{}2*c\^{}2+396*b*w\^{}2*c\^{}3\relax
-288*b*w*c\^{}3)*b/((w-1)(44*b\^{}2*c\^{}2+36\newline
*b*c\^{}2+22*b\^{}2*c\^{}4+68*b\^{}2*c-134*b\^{}2*c\^{}3\relax
+150*b\^{}3*c\^{}4+160*b\^{}3*c-166*b\^{}4\newline
*c\^{}4+490*b\^{}4*c\^{}3-332*b\^{}4*c\^{}2-116*b\^{}4*c\relax
+40*b\^{}3+136*b\^{}4+6*c\^{}3-18*c\^{}4*b\newline
-300*b\^{}3*c\^{}2+17*c\^{}5*b\^{}2-40*c\^{}5*b\^{}3\relax
+17*c\^{}6*b\^{}4-29*c\^{}5*b\^{}4-5*c\^{}6*b\^{}3+3\newline
96*b\^{}5*c\^{}2-432*b\^{}5*c+198*b\^{}6*c\^{}4\relax
-378*b\^{}6*c\^{}3+396*b\^{}6*c\^{}2-216*b\^{}6*c+1\newline
08*b\^{}5*c\^{}5+144*b\^{}5+48*b\^{}6+6*b\^{}6*c\^{}6\relax
-54*b\^{}6*c\^{}5-18*b\^{}5*c\^{}6-198*b\^{}5*c\^{}\newline
4)(b\^{}2*c\^{}4-6*b\^{}2*c\^{}3+13*b\^{}2*c\^{}2-12*b\^{}2*c\relax
+4*b\^{}2+c\^{}2));}
\medskip\noindent
{\tt 
d3=-1/9*(-9882*w\^{}2*b\^{}5*c\^{}6-81*w\^{}2*c\^{}9*b\^{}3\relax
-13392*b\^{}5*w\^{}2*c+243*w\^{}2*c\^{}\newline
8*b\^{}3+2232*b\^{}2*c\^{}2-4752*b\^{}2*c\^{}4\relax
-2016*b\^{}2*c\^{}3+16236*b\^{}3*c\^{}4+1044*b*c\^{}\newline
3+1584*b\^{}3*c+17064*b\^{}4*c\^{}4+576*b\^{}4*c\^{}3\relax
-6408*b\^{}4*c\^{}2+864*b\^{}4*c+864*b\^{}\newline
4-1368*c\^{}4*b-15408*b\^{}3*c\^{}3-279*c\^{}8*b\^{}3\relax
+2232*b\^{}3*c\^{}2-2376*c\^{}6*b\^{}2+716\newline
4*c\^{}5*b\^{}2+3852*c\^{}7*b\^{}3-801*c\^{}8*b\^{}4\relax
+144*c\^{}7*b\^{}4+8532*c\^{}6*b\^{}4-504*c\^{}7*\newline
b\^{}2-20916*c\^{}5*b\^{}4-8118*c\^{}6*b\^{}3+684*c\^{}6*b\relax
+7920*b\^{}6*c\^{}8-29412*b\^{}5*c\^{}3+\newline
28728*b\^{}5*c\^{}2-14832*b\^{}5*c+144900*b\^{}6*c\^{}4\relax
-118080*b\^{}6*c\^{}3+63360*b\^{}6*c\^{}\newline
2-20160*b\^{}6*c-3591*b\^{}5*c\^{}8+7353*b\^{}5*c\^{}7\relax
+180*c\^{}4-180*c\^{}5+90*c\^{}6+54*w\^{}\newline
2*c\^{}10*b\^{}6-29520*b\^{}6*c\^{}7+3168*b\^{}5\relax
+2880*b\^{}6+72450*b\^{}6*c\^{}6-122580*b\^{}6*\newline
c\^{}5-7758*b\^{}5*c\^{}6-70848*c\^{}3*b\^{}6*w\^{}2\relax
-3645*w\^{}2*b\^{}5*c\^{}8+137*w\^{}2*c\^{}8*b\^{}2+\newline
15516*b\^{}5*c\^{}4+927*c\^{}9*b\^{}5-1260*c\^{}9*b\^{}6\relax
-99*c\^{}9*b\^{}3+54*c\^{}9*b\^{}4+279*c\^{}8\newline
*b\^{}2+90*c\^{}10*b\^{}6+27*c\^{}10*b\^{}4\relax
+4752*w\^{}2*b\^{}6*c\^{}8+837*w\^{}2*c\^{}9*b\^{}5-73548*\newline
w\^{}2*b\^{}6*c\^{}5+108*w\^{}2*c\^{}4-756*w\^{}2*c\^{}9*b\^{}6\relax
-1008*b\^{}2*w\^{}2*c\^{}4+29160*c\^{}2*b\newline
\^{}5*w\^{}2+2888*w\^{}2*c\^{}7*b\^{}4\relax
-17712*w\^{}2*b\^{}6*c\^{}7-6916*w\^{}2*c\^{}6*b\^{}4+1988*w\^{}2*\newline
c\^{}5*b\^{}2-3896*b\^{}4*w\^{}2*c\^{}2+8451*w\^{}2*b\^{}5*c\^{}7\relax
-12096*c*b\^{}6*w\^{}2-487*w\^{}2*c\^{}\newline
8*b\^{}4+43470*w\^{}2*b\^{}6*c\^{}6+1728*b\^{}6*w\^{}2\relax
+54*w\^{}2*c\^{}6-1344*b\^{}2*w\^{}2*c\^{}3+29*\newline
w\^{}2*c\^{}10*b\^{}4+86940*w\^{}2*b\^{}6*c\^{}4\relax
+324*w\^{}2*c\^{}7*b\^{}3-78*w\^{}2*c\^{}9*b\^{}4+2916*c\newline
\^{}4*b\^{}3*w\^{}2-1248*b\^{}4*w\^{}2*c\relax
-1944*b\^{}3*w\^{}2*c\^{}2+2592*b\^{}5*w\^{}2\relax
-81*b\^{}5*w\^{}2*c\newline
\^{}10-108*w\^{}2*c\^{}5-648*b*w\^{}2*c\^{}4\relax
-99*c\^{}10*b\^{}5+38016*c\^{}2*b\^{}6*w\^{}2-1296*b\^{}3\newline
*w\^{}2*c\^{}3+324*w\^{}2*c\^{}6*b-504*c\^{}6*b\^{}2*w\^{}2\relax
-336*c\^{}7*b\^{}2*w\^{}2-1458*c\^{}6*b\^{}3*\newline
w\^{}2-33804*c\^{}3*b\^{}5*w\^{}2+19764*c\^{}4*b\^{}5*w\^{}2\relax
-13832*c\^{}4*b\^{}4*w\^{}2+11060*c\^{}5*\newline
b\^{}4*w\^{}2+11552*c\^{}3*b\^{}4*w\^{}2+928*b\^{}4*w\^{}2\relax
-135*b*w\^{}2*c\^{}7-261*c\^{}7*b+1296*b\newline
\^{}3*w\^{}2*c+1096*b\^{}2*w\^{}2*c\^{}2\relax
+540*b*w\^{}2*c\^{}3)*b/((1+w)(w-1)(44*b\^{}2*c\^{}2+36\newline
*b*c\^{}2+22*b\^{}2*c\^{}4+68*b\^{}2*c-134*b\^{}2*c\^{}3\relax
+150*b\^{}3*c\^{}4+160*b\^{}3*c-166*b\^{}4\newline
*c\^{}4+490*b\^{}4*c\^{}3-332*b\^{}4*c\^{}2-116*b\^{}4*c\relax
+40*b\^{}3+136*b\^{}4+6*c\^{}3-18*c\^{}4*b\newline
-300*b\^{}3*c\^{}2+17*c\^{}5*b\^{}2-40*c\^{}5*b\^{}3\relax
+17*c\^{}6*b\^{}4-29*c\^{}5*b\^{}4-5*c\^{}6*b\^{}3+3\newline
96*b\^{}5*c\^{}2-432*b\^{}5*c+198*b\^{}6*c\^{}4\relax
-378*b\^{}6*c\^{}3+396*b\^{}6*c\^{}2-216*b\^{}6*c+1\newline
08*b\^{}5*c\^{}5+144*b\^{}5+48*b\^{}6+6*b\^{}6*c\^{}6\relax
-54*b\^{}6*c\^{}5-18*b\^{}5*c\^{}6-198*b\^{}5*c\^{}\newline
4)(b\^{}2*c\^{}4-6*b\^{}2*c\^{}3+13*b\^{}2*c\^{}2-12*b\^{}2*c\relax
+4*b\^{}2+c\^{}2)).}\par

\enddocument
\end